\documentclass[12pt]{article}%
\usepackage{amsmath}
\usepackage{amsfonts}
\usepackage{amssymb}
\usepackage{graphicx}%
\begin{document}
\vspace*{.5in}

\centerline{{\large \textbf{Order-Reducing Form Symmetries and Semiconjugate} }}
\vspace{0.35ex}
\centerline{{\large \textbf{Factorizations of Difference Equations} }}
\vspace{2ex}

\centerline{H. SEDAGHAT} \vspace{3ex}

\begin{quote}
\noindent{\small \textbf{Abstract.} The scalar difference equation
$x_{n+1}=f_{n}(x_{n},x_{n-1},\cdots,x_{n-k})$ may exhibit symmetries in its
form that allow for reduction of order through substitution or a change of
variables. Such form symmetries can be defined generally using the
semiconjugate relation on a group which yields a reduction of order through
the semiconjugate factorization of the difference equation of order
}${\small k+1}$ {\small into equations of lesser orders. Different classes of
equations are considered including separable equations and homogeneous
equations of degree 1. Applications include giving a complete factorization of
the linear non-homogeneous difference equation of order }${\small k+1}$
{\small into a system of }${\small k+1}${\small \ first order linear
non-homogeneous equations in which the coefficients are the eigenvalues of the
higher order equation. Form symmetries are also used to explain the
complicated multistable behavior of a separable, second order exponential
equation.}
\end{quote}

\medskip

\noindent\noindent{\small \textbf{Keywords.} Form symmetry, order reduction,
semiconjugate, groups, difference equations, linear non-homogeneous,
separable, homogeneous of degree 1, multistability}

\vspace{2ex}

\section{Introduction}

Certain difference equations have symmetries in their expressions that allow a
reduction of their orders through substitutions of new variables. For
instance, consider the second order, scalar difference equation%
\begin{equation}
x_{n+1}=x_{n}+g_{n}(x_{n}-x_{n-1}) \label{I1}%
\end{equation}
where $g_{n}$ is a real function for each integer $n$. This equation has a
symmetry in its form that is easy to identify when (\ref{I1}) is re-written as%
\begin{equation}
x_{n+1}-x_{n}=g_{n}(x_{n}-x_{n-1}). \label{I2}%
\end{equation}

Now, setting $t_{n}=x_{n}-x_{n-1}$ changes Eq.(\ref{I2}) to the first order
equation%
\begin{equation}
t_{n+1}=g_{n}(t_{n}). \label{I3}%
\end{equation}

The expression $x_{n}-x_{n-1}$ is an example of what we may call a
\textit{form symmetry}. Substituting a new variable $t_{n}$ for this form
symmetry in (\ref{I2}) gave the lower order equation (\ref{I3}). The form
symmetry also establishes a link between the second order equation and the
first order one, in the sense that information about each solution $\{t_{n}\}$
of (\ref{I3}) can then be translated into information about the corresponding
solution of (\ref{I1}) using the equation%
\begin{equation}
x_{n}=x_{n-1}+t_{n}=x_{0}+\sum_{k=1}^{n}t_{k} \label{I3a}%
\end{equation}
where $x_{0}$ is an initial value for (\ref{I1}). Along similar lines, the
non-homogeneous linear difference equation
\begin{equation}
x_{n+1}+px_{n}+qx_{n-1}=\alpha_{n}\quad\text{with }1+p+q=0 \label{lin2}%
\end{equation}
has at least two form symmetries. First, setting $p=-1-q$ and rearranging
terms in (\ref{lin2}) reveals the form symmetry $x_{n}-x_{n-1}$ and the
corresponding order-reducing substitution:%
\begin{equation}
x_{n+1}-x_{n}=\alpha_{n}+q(x_{n}-x_{n-1})\Rightarrow t_{n+1}=\alpha_{n}%
+qt_{n}. \label{lin2a}%
\end{equation}

Further,%
\begin{equation}
x_{n+1}-qx_{n}=\alpha_{n}-(p+q)x_{n}-qx_{n-1}=\alpha_{n}+x_{n}-qx_{n-1}%
\Rightarrow t_{n+1}=\alpha_{n}+t_{n}. \label{lin2b}%
\end{equation}

Thus $x_{n}-qx_{n-1}$ is also a form symmetry of (\ref{lin2}). Note that the
coefficients $q$ and $1$ of $t_{n}$ in (\ref{lin2a}) and (\ref{lin2b}),
respectively, are both eigenvalues of the homogeneous part of (\ref{lin2}),
i.e., roots of the characteristic polynomial $z^{2}+pz+q$ when $p+q=-1$. Later
in this article we show that this relationship holds for all linear difference equations.

These and many other types of known form symmetries of difference equations of
type%
\begin{equation}
x_{n+1}=f_{n}(x_{n},x_{n-1},\cdots,x_{n-k}). \label{kth}%
\end{equation}
can be defined in terms of semiconjugacy; in \cite{hsbk} there is a basic
discussion of this topic for real functions but we give a more general
definition in this article. The idea of reducing order via semiconjugacy is
basically simple; we find functions that are semiconjugate to the given
functions $f_{n}$ but which have fewer variables than $k+1,$ i.e., the order
of (\ref{kth}). Then the semiconjugate relation simultaneously determines both
a form symmetry and a factorization of Eq.(\ref{kth}). This factorization of
difference equations, which is also a formal representation of the
substitution process discussed above, yields a pair of lower order equations.
This pair is made up of a factor equation such as (\ref{I3}) and an associated
equation such as (\ref{I3a}) that is derived from the form symmetry and
relates the factor equation to the original one. The orders of the factor
equation and the associated one always add up to the order of the original
scalar equation.

The aim of this article is to formalize, within the framework of
semiconjugacy, the concept of form symmetry and its use in reduction of order
by substitutions. In addition to unifying various ad hoc techniques, this
approach also gives rise to new methods for analyzing higher order difference
equations. Some of these new methods are described in this article along with
examples and applications.

Symmetries of a different kind have already been used to study higher order
difference equations or systems of first order ones. A well-known approach
involves adaptation of the Lie symmetry concept from differential equations to
difference equations; see, e.g., \cite{hydon} for a general discussion of the
discrete case covering various topics such as reduction of order,
integrability and finding explicit solutions. Also see \cite{bsq}, \cite{ltw},
\cite{maeda} for additional ideas and techniques. The main difference between
the concepts of form symmetry and Lie symmetry may be summed up as follows:
Form symmetries are sought in the difference equation itself whereas Lie
symmetries exist in the solutions of the equation. The existence of either
type of symmetry can yield valuable information about the dynamics and the
solutions of the difference equation with a variety of applications such as
reduction of order.

\section{Semiconjugate forms}

The material in this section substantially extends the notions in \cite{hsbk}
to the more general group context. For related results and background material
on difference equations see, e.g., \cite{ag}, \cite{bc}, \cite{ce}, \cite{el},
\cite{kl}.

In this article, $G$ denotes a non-trivial group. The group structure provides
a suitable framework for our results. However, in most applications $G$ turns
out to be a substructure of\ a more complex object such as a vector space, a
ring or an algebra possessing a compatible or natural metric topology. In
typical studies involving difference equations and discrete dynamical systems,
$G$ is a group of real or complex numbers. We may define each mapping $f_{n}$
on the ambient sturcture as long as the following invariance condition holds%
\begin{equation}
f_{n}(G^{k+1})\subset G\text{ \ for all }n. \label{invar}%
\end{equation}

If (\ref{invar}) holds then for each set of initial values $x_{0}%
,x_{-1},\ldots,x_{-k}$ in $G$ Eq.(\ref{kth}) recursively generates a
\textit{solution} or \textit{orbit} $\{x_{n}\}_{n=-1}^{\infty}$ in $G.$ Before
defing form symmetries for (\ref{kth}), it is necessary to discuss some
general defintions involving systems.

Let $1\leq m\leq k.$ We say that a self map $F_{n}=[f_{1,n},\ldots,f_{k+1,n}]$
of $G^{k+1}$ is \textit{semiconjugate} to a self map $\Phi_{n}$ of $G^{m}$ if
there is a function $H:G^{k+1}\rightarrow G^{m}$ such that for every $n,$%
\begin{equation}
H\circ F_{n}=\Phi_{n}\circ H. \label{sc0}%
\end{equation}

Each mapping $\Phi_{n}$ is called a semiconjugate \textit{factor} of the
corresponding $F_{n}$. Suppose that%
\begin{align*}
H(u_{0},\ldots,u_{k})  &  =[h_{1}(u_{0},\ldots,u_{k}),\ldots,h_{m}%
(u_{0},\ldots,u_{k})]\\
\Phi_{n}(t_{1},\ldots,t_{m})  &  =[\phi_{1,n}(t_{1},\ldots,t_{m}),\ldots
,\phi_{m,n}(t_{1},\ldots,t_{m})]
\end{align*}
where $h_{j}:G^{k+1}\rightarrow G$ and $\phi_{j,n}:G^{m}\rightarrow G$ are the
corresponding component functions. Then identity (\ref{sc0}) is equivalent to
the system%
\begin{align}
h_{j}(f_{1,n}(u_{0},\ldots,u_{k}),\ldots,f_{k+1,n}(u_{0},\ldots,u_{k}))  &
=\phi_{j,n}(h_{1}(u_{0},\ldots,u_{k}),\ldots,h_{m}(u_{0},\ldots,u_{k}%
))\nonumber\\
j  &  =1,2,\ldots,m. \label{scj}%
\end{align}

If the functions $f_{j,n}$ are given then (\ref{scj}) is a system of
functional equations whose solutions $h_{j},\phi_{j,n}$ give the maps $H$ and
$\Phi_{n}.$ The functions $\Phi_{n}$ on $G^{m}$ define a system with lower
dimension than that defined by the functions $F_{n}$ on $G^{k+1}.$ For a given
solution $\{X_{n}\}$ of the equation%
\begin{equation}
X_{n+1}=F_{n}(X_{n}),\quad X_{0}\in G^{k+1} \label{Fk}%
\end{equation}
let $Y_{n}=H(X_{n})$ for $n=0,1,2,\ldots$ Then
\[
Y_{n+1}=H(X_{n+1})=H(F_{n}(X_{n}))=\Phi_{n}(H(X_{n}))=\Phi_{n}(Y_{n})
\]
so that $\{Y_{n}\}$ satisfies the lower order equation
\begin{equation}
Y_{n+1}=\Phi_{n}(Y_{n}),\quad Y_{0}=H(X_{0})\in G^{m}. \label{phim}%
\end{equation}

The relationship between the solutions of (\ref{Fk}) and those of (\ref{phim})
is not generally straightforward; however, information about the solutions of
(\ref{phim}) can shed light on the dynamics of (\ref{Fk}). The case where
$G=\mathbb{R}$, $m=1$ and $F_{n}=F$ is time-independent (or autonomous) is of
some interest because in this case Eq.(\ref{phim}) is a first order difference
equation on the real numbers and as such its dynamics are much better
understood than that of the higher dimensional Eq.(\ref{Fk}). This case is
discussed in detail in \cite{hsbk}.

\section{Order-reducing form symmetries}

For the scalar difference equation (\ref{kth}) that is of interest in this
article, each $F_{n}$ is the associated vector map (or unfolding) of the
function $f_{n}$ in Eq.(\ref{kth}), i.e.,%
\[
F_{n}(u_{0},\ldots,u_{k})=[f_{n}(u_{0},\ldots,u_{k}),u_{0},\ldots,u_{k-1}].
\]

Even if each such $F_{n}$ is semiconjugate to an $m$-dimensional map $\Phi
_{n}$ as in (\ref{sc0}), the preceding discussion only gives the system
(\ref{phim}) in which the maps $\Phi_{n}$ are not necessarily of scalar type
similar to $F_{n}.$ While this may be unavoidable in some cases, adding a few
reasonable restrictions can ensure that each $\Phi_{n}$ is also of scalar
type. To this end, define%
\begin{equation}
h_{1}(u_{0},\ldots,u_{k})=u_{0}\ast h(u_{1},\ldots,u_{k}) \label{H}%
\end{equation}
where $h:G^{k}\rightarrow G$ is a function to be determined and $\ast$ denotes
the group operation. This restriction on $H$ makes sense for Eq.(\ref{kth}),
which is of recursive type; i.e., $x_{n+1}$ given explicitly by functions
$f_{n}$. With these restrictions on $H$ and $F_{n}$ the first equation in
(\ref{scj}) is given by%
\begin{equation}
f_{n}(u_{0},\ldots,u_{k})\ast h(u_{0},\ldots,u_{k-1})=g_{n}(u_{0}\ast
h(u_{1},\ldots,u_{k}),\ldots,h_{m}(u_{0},\ldots,u_{k})) \label{scm}%
\end{equation}
where for notational convenience we have set%
\[
g_{n}\doteq\phi_{1,n}:G^{m}\rightarrow G.
\]

Eq.(\ref{scm}) is a functional equation in which the functions $h,h_{j},g_{n}$
may be determined in terms of the given functions $f_{n}$. Our aim is
ultimately to extract a scalar equation of order $m$ such as%
\begin{equation}
t_{n+1}=g_{n}(t_{n},\ldots,t_{n-m+1}) \label{gnt}%
\end{equation}
from (\ref{scm}) in such a way that the maps $\Phi_{n}$ will be of scalar
type. The basic framework is already in place; let $\{x_{n}\}$ be a solution
of Eq.(\ref{kth}) and define%
\[
t_{n}=x_{n}\ast h(x_{n-1},\ldots,x_{n-k}).
\]

Then the left hand side of\ (\ref{scm}) is%
\[
x_{n+1}\ast h(x_{n},\ldots,x_{n-k+1})=t_{n+1}.
\]
which gives the initial part of the difference equation (\ref{gnt}). In order
that the right hand side of (\ref{scm}) coincide with that in (\ref{gnt}) it
is necessary to define%
\begin{equation}
h_{j}(x_{n},\ldots,x_{n-k})=t_{n-j+1}=x_{n-j+1}\ast h(x_{n-j},\ldots
,x_{n-k-j+1}),\quad j=2,\ldots,m. \label{hgj}%
\end{equation}

Since the left hand side of (\ref{hgj}) does not depend on terms
$x_{n-k-1},\ldots,x_{n-k-j+1}$ it follows that the function $h_{j}$ must be
constant in its last few coordinates. Since $h$ deos not depend on $j$ the
number of constant coordinates is found from the last function $h_{m}.$
Specifically, we have%
\begin{equation}
h_{m}(x_{n},\ldots,x_{n-k})=x_{n-m+1}\ast h(\underset{k-m+1\text{ variables}%
}{\underbrace{x_{n-m},\ldots,x_{n-k}}},\underset{m-1\text{ terms with }h\text{
constant}}{\underbrace{x_{n-k-1},\ldots,x_{n-k-m+1}}}) \label{consis}%
\end{equation}

The preceding condition leads to the necessary restrictions on $h$ and every
$h_{j}$ for a consistent derivation of (\ref{gnt}) from (\ref{scm}), so
(\ref{consis}) is a consistency condition. Now from (\ref{hgj}) and
(\ref{consis}) we obtain for $(u_{0},\ldots,u_{k})\in G^{k+1}$
\begin{equation}
h_{j}(u_{0},\ldots,u_{k})=u_{j-1}\ast h(u_{j},u_{j+1}\ldots,u_{j+k-m}),\quad
j=1,\ldots,m. \label{hjd}%
\end{equation}

We refer to $H=[h_{1},\ldots,h_{m}]$ as a\textit{ form symmetry} for
Eq.(\ref{kth}) if the components $h_{j}$ are defined by (\ref{hjd}). Since the
range of $H$ has a lower dimension than its domain, we say that $H$ is an
\textit{order-reducing} form symmetry.

Using the forms in (\ref{hjd}) for $h_{j}$ in (\ref{scm}) for every solution
$\{x_{n}\}$ of Eq.(\ref{kth}) we obtain the following pair of equations from
(\ref{scm}), the first of which is just (\ref{gnt}):
\begin{subequations}
\label{scf}%
\begin{align}
t_{n+1}  &  =g_{n}(t_{n},\ldots,t_{n-m+1}),\label{fac}\\
x_{n+1}  &  =t_{n+1}\ast h(x_{n},\ldots,x_{n-k+m})^{-1}. \label{cofac}%
\end{align}

The power $-1$ represents group inversion in $G.$ The first equation
(\ref{fac}) may be called a \textit{factor} of Eq.(\ref{kth}) since it is
distilled from the semiconjugate factor $\Phi_{n}.$ The second equation
(\ref{cofac}) that links the factor to the original equation may be called a
\textit{cofactor} of Eq.(\ref{kth}). We call the system of equations
(\ref{scf}) a \textit{semiconjugate (SC)} \textit{factorization} of
Eq.(\ref{kth}).

Note that if $\{t_{n}\}$ is a given solution of (\ref{fac}) then using this
sequence in (\ref{cofac}) produces a solution $\{x_{n}\}$ of (\ref{kth}).
Conversely, if $\{x_{n}\}$ is a solution of (\ref{kth}) then the sequence
$t_{n}=x_{n}\ast h(x_{n-1},\ldots,x_{n-k+m-1})$ is a solution of (\ref{fac})
with initial values
\end{subequations}
\[
t_{-j}=x_{-j}\ast h(x_{-j-1},\ldots,x_{-j-k+m-1}),\quad j=0,\ldots,m-1.
\]

Since solutions of the pair of equations (\ref{fac}) and (\ref{cofac})
coincide with the solutions of the scalar equation (\ref{kth}), we say that
the pair (\ref{fac}) and (\ref{cofac}) is equivalent to (\ref{kth}). The
following summarizes the preceding discussions.

\medskip

\noindent\textbf{Theorem 1}. \textit{Let }$k\geq1$\textit{, }$1\leq m\leq
k$\textit{ and suppose that there are functions }$h:G^{k-m+1}\rightarrow
G$\textit{ and }$g_{n}:G^{m}\rightarrow G$\textit{ that satisfy equations
(\ref{scm}) and (\ref{hjd}). Then with the order-reducing form symmetry }%
\[
H(u_{0},\ldots,u_{k})=[u_{0}\ast h(u_{1},\ldots,u_{k+1-m}),\ldots,u_{m-1}\ast
h(u_{m},\ldots,u_{k})]
\]
\textit{Eq.(\ref{kth}) is equivalent to the SC factorization consisting of the
pair of equations (\ref{fac}) and (\ref{cofac}) whose orders }$m$\textit{ and
}$k+1-m$\textit{ respectively, add up to the order of (\ref{kth}).}

\medskip

In this setting we say that the SC factorization (\ref{scf}) gives a
\textit{type-(}$m,k+1-m$\textit{) order reduction} for Eq.(\ref{kth}), or that
(\ref{kth}) is a type-($m,k+1-m$) equation. A second order difference equation
($k=1$) can have only the order-reduction type (1,1) into two first order
equations although the factor and cofactor equations are not uniquely defined.
In general, a higher order difference equation may have more than one SC
factorization. A third order equation can have two order-reduction types,
namely, (1,2) and (2,1). Of the $k$ possible order reduction types for an
equation of order $k+1$ the two extreme ones, namely, $(1,k)$ and $(k,1)$ have
the extra appeal of having an equation of order 1 as either a factor or a
cofactor. In the next two sections we discuss classes of higher order
difference equations having one of these order-reduction types.

We note that the SC factorization of Theorem 1 does not require the
determination of $\phi_{j,n}$ for $j\geq2.$ For completeness, we close this
section by showing that each coordinate function $\phi_{j,n}$ projects into
coordinate $j-1$ for $j>1,$ thus showing that $\Phi_{n}$ is of scalar type,
i.e., it is the unfolding of Eq.(\ref{gnt}) in the same sense that $F_{n}$
unfolds (\ref{kth}). If the maps $h_{j}$ are given by (\ref{hjd}) then for
$j\geq2$ (\ref{scj}) gives%
\begin{align*}
\phi_{j,n}(h_{1}(u_{0},\ldots,u_{k}),\ldots,h_{m}(u_{0},\ldots,u_{k})) &
=h_{j}(f_{n}(u_{0},\ldots,u_{k}),u_{0},\ldots,u_{k-1})\\
&  =u_{j-2}\ast h(u_{j-1},u_{j}\ldots,u_{j+k-m-1})\\
&  =h_{j-1}(u_{0},\ldots,u_{k}).
\end{align*}

Therefore, for each $n$ and for every $(t_{1},\ldots,t_{m})\in H(G^{k+1})$ we
have
\[
\Phi_{n}(t_{1},\ldots,t_{m})=[g_{n}(t_{1},\ldots,t_{m}),t_{1},\ldots,t_{m-1}]
\]
i.e., $\Phi_{n}|_{H(G^{k+1})}$ is of scalar type. Further, if $H$ is defined
component-wise by (\ref{hjd}) then $H(G^{k+1})=G^{m};$ i.e., $H$ is onto
$G^{m}$ so that $\Phi_{n}$ is of scalar type. To prove the onto claim, we pick
arbitrary $[t_{1},\ldots,t_{m}]\in G^{m}$ and set $u_{m-1}=t_{m}\ast
h(u_{m},u_{m+1}\ldots,u_{k})^{-1}$ where $u_{m}=u_{m+1}=\ldots u_{k}=1$ (the
group identity).\ Then
\begin{align*}
t_{m} &  =u_{m-1}\ast h(1,1\ldots,1)\\
&  =u_{m-1}\ast h(u_{m},u_{m+1}\ldots,u_{k})\\
&  =h_{m}(u_{0},\ldots,u_{k})\\
&  =h_{m}(u_{0},\ldots,u_{m-2},t_{m}\ast h(1,1\ldots,1)^{-1},1\ldots,1).
\end{align*}
for any choice of $u_{0},\ldots,u_{m-2}\in G.$ Similarly, define
$u_{m-2}=t_{m-1}\ast h(u_{m-1},u_{m}\ldots,u_{k-1})^{-1}$ so as to get
\begin{align*}
t_{m-1} &  =u_{m-2}\ast h(u_{m-1},u_{m}\ldots,u_{k-1})\\
&  =h_{m-1}(u_{0},\ldots,u_{k})\\
&  =h_{m-1}(u_{0},\ldots,u_{m-3},t_{m-1}\ast h(u_{m-1},1\ldots,1)^{-1}%
,u_{m-1},1\ldots,1)
\end{align*}
for any choice of $u_{0},\ldots,u_{m-3}\in G.$ Continuing in this way,
induction leads to selection of $u_{m-1},\ldots,u_{0}$ such that%
\[
t_{j}=h_{j}(u_{0},\ldots,u_{m-1},1,\ldots,1),\quad j=1,\ldots,m
\]
and it is proved that $H$ is onto $G^{m}.$

\section{HD1 and other type-$(k,1)$ factorizations}

If $m=k$ then the function $h:G\rightarrow G$ in (\ref{hjd}) is of one
variable and we obtain a type-$(k,1)$ order reduction with form symmetry%
\begin{equation}
H(u_{0},\ldots,u_{k})=[u_{0}\ast h(u_{1}),u_{1}\ast h(u_{2})\ldots,u_{k-1}\ast
h(u_{k})] \label{km}%
\end{equation}
and SC factorization%
\begin{align*}
t_{n+1}  &  =g_{n}(t_{n},t_{n-1},\ldots,t_{n-k+1})\\
x_{n+1}  &  =t_{n+1}\ast h(x_{n})^{-1}%
\end{align*}
where the functions $g_{n}:G^{k}\rightarrow G$ are determined by the given
functions $f_{n}$ in (\ref{kth}) as in the previous section.

The simplest example of a non-constant $h$ in this setting is the identity
function $h(u)=u$ for all $u\in G.$ An example of a type-$(k,1)$ difference
equation having this type of form symmetry over $(0,\infty)$ under ordinary
multiplication is the rational equation%
\begin{equation}
x_{n+1}=\frac{ax_{n-1}}{x_{n}x_{n-1}+b},\quad a,b>0. \label{k1rat}%
\end{equation}

The term $x_{n}x_{n-1}$ in the denominator suggests multiplying (\ref{k1rat})
by $x_{n}$ on both sides and substituting
\[
t_{n}=x_{n}x_{n-1}=x_{n}h(x_{n-1})
\]
to get the SC factorization%
\begin{align}
t_{n+1}  &  =\frac{at_{n}}{t_{n}+b}=g(t_{n}),\quad t_{0}=x_{0}x_{-1}%
\label{ric}\\
x_{n+1}  &  =\frac{t_{n+1}}{x_{n}}=\frac{t_{n+1}}{h(x_{n})}.\nonumber
\end{align}

Further, Eq.(\ref{ric}) can be made linear by the change of variables
$s_{n}=1/t_{n}.$ For an exhaustive treatment of (\ref{k1rat}) based on these
ideas, see \cite{ord23}. Another form symmetry of type (\ref{km}) that is
defined on $\mathbb{C}$ or $\mathbb{R}$ is based on $h(u)=cu$ where $c$ is a
fixed, nonzero complex or real number. This type of form symmetry (with real
$c$) has been used in e.g., \cite{dkmos2} and \cite{kyoto}.

In the case where $h(u)=u^{-1}$ is based on group inversion, it is possible to
identify the class of functions $f_{n}$ that have the form symmetry
(\ref{km}). Equation (\ref{kth}) is said to be \textit{homogeneous of degree 1
(HD1)} if for every $n=1,2,3,\ldots$ the functions $f_{n}$ are homogeneous of
degree 1 relative to the group $G,$ i.e.,%
\[
f_{n}(u_{0}\ast t,\ldots,u_{k}\ast t)=f_{n}(u_{0},\ldots,u_{k})\ast t\text{
for all }t,u_{i}\in G,\ i=0,\ldots,k,\ n\geq1.
\]

If $G$ is non-commutative then this definition gives a \textquotedblleft right
version\textquotedblright\ of the HD1 property; a \textquotedblleft left
version\textquotedblright\ can be defined analogously. We note that the two
equations (\ref{I1}) and (\ref{lin2}) in the Introduction are HD1 relative to
the additive group of real numbers. For comments on homogeneous functions and
their abundance on groups we refer to \cite{hd1}; though stated for functions
of two variables, the results in \cite{hd1} easily extend to any number of
variables. The following result shows that the HD1 property characterizes the
inversion-based form symmetry and yields a type-$(k,1)$ order-reduction in
every case.

\medskip

\noindent\textbf{Theorem 2}. \textit{Eq.(\ref{kth}) has the inversion-based
form symmetry }%
\begin{equation}
H(u_{0},\ldots,u_{k})=[u_{0}\ast u_{1}^{-1},\ldots,u_{k-1}\ast u_{k}%
^{-1}],\quad h(t)=t^{-1}\label{infs}%
\end{equation}
\textit{if and only if }$f_{n}$\textit{ is HD1 relative to }$G$\textit{ for
all }$n$. \textit{In this case, (\ref{kth}) has a type-(}$k,1$\textit{)
order-reduction with the SC factorization}
\begin{subequations}
\label{1a}%
\begin{align}
t_{n+1} &  =f_{n}(1,t_{n}^{-1},(t_{n}\ast t_{n-1})^{-1},\ldots,(t_{n}\ast
t_{n-1}\ast\cdots\ast t_{n-k+1})^{-1})\label{1aa}\\
x_{n+1} &  =t_{n+1}\ast x_{n}.\label{1ab}%
\end{align}

\textit{Note that the factor difference equation (\ref{1aa}) has order }%
$k$\textit{ and its cofactor (\ref{1ab}) is linear non-autonomous of order one
in }$x_{n}$\textit{.}

\textbf{Proof.} First, assume that (\ref{kth}) has the form symmetry
(\ref{infs}) that satisfies Eq.(\ref{scm}) for given functions $g_{n},$ i.e.,%
\end{subequations}
\begin{equation}
f_{n}(u_{0},\ldots,u_{k})\ast u_{0}^{-1}=g_{n}(u_{0}\ast u_{1}^{-1}%
,\ldots,u_{k-1}\ast u_{k}^{-1}).\label{hd1a}%
\end{equation}

Let $t\in G$ be arbitrary. Then for all $n$ (\ref{hd1a}) implies%
\begin{align*}
f_{n}(u_{0}\ast t,\ldots,u_{k}\ast t)  & =g_{n}((u_{0}\ast t)\ast(u_{1}\ast
t)^{-1},\ldots,(u_{k-1}\ast t)\ast(u_{k}\ast t)^{-1})\ast(u_{0}\ast t)\\
& =[g_{n}(u_{0}\ast u_{1}^{-1},\ldots,u_{k-1}\ast u_{k}^{-1})\ast u_{0}]\ast
t\\
& =f_{n}(u_{0},\ldots,u_{k})\ast t.
\end{align*}

It follows that $f_{n}$ is HD1 relative to $G$ for all $n$ and the
first part of the theorem is proved. The converse is proved in a
straightforward fashion; see \cite{hd1a}.

\medskip

\noindent\textbf{Remarks. }

1. Equation (\ref{1ab}) can be solved explicitly in terms of a solution
$\{t_{n}\}$ of (\ref{1aa}) as follows:
\begin{equation}
x_{n}=%
%TCIMACRO{\tprod _{i=0}^{n-1}}%
%BeginExpansion
{\textstyle\prod_{i=0}^{n-1}}
%EndExpansion
t_{n-i}\ast x_{0}\quad n=1,2,3,\ldots\label{3}%
\end{equation}
where the multiplicative notation is used for iterations of the group
opreation $\ast.$ In additive (and commutative) notation, (\ref{3}) takes the
form%
\begin{equation}
x_{n}=x_{0}+%
%TCIMACRO{\tsum _{i=1}^{n}}%
%BeginExpansion
{\textstyle\sum_{i=1}^{n}}
%EndExpansion
t_{i}. \label{3a}%
\end{equation}

2. We can quickly construct Eq.(\ref{1aa}) directly from (\ref{kth}) in the
HD1 case by making the substitutions%
\begin{equation}
1\rightarrow x_{n},\quad(t_{n}t_{n-1}\cdots t_{n-i+1})^{-1}\rightarrow
x_{n-i}\ \ \text{for }i=1,2,\ldots,k. \label{subs}%
\end{equation}

Recall that 1 represents the group identity in multiplicative notation. In
additive notation (\ref{subs}) takes the form%
\begin{equation}
0\rightarrow x_{n},\quad-t_{n}-t_{n-1}\cdots-t_{n-i+1}\rightarrow
x_{n-i}\ \ \text{for }i=1,2,\ldots,k. \label{subsum}%
\end{equation}

\medskip

Previous studies involving HD1 equations implicitly use the idea behind
Theorem 2 above to reduce second order equations to first order ones; see e.g.
\cite{dkmos1}, \cite{ks}, \cite{ord23}. Examples 1-3 next illustrate Theorem 2
and some associated concepts.

\medskip

\noindent\textbf{Example 1}. Consider the rational delay difference equation%
\begin{equation}
x_{n+1}=x_{n}\left(  \frac{a_{n}x_{n-k+1}}{x_{n-k}}+b_{n}\right)  ,\quad
x_{0},x_{-1},\ldots,x_{-k}>0\label{rk}%
\end{equation}
where $\{a_{n}\},\{b_{n}\}$ are sequences of positive real numbers. This
equation is HD1 relative to the group $(0,\infty)$ under ordinary
multiplication. Thus Theorem 2 and (\ref{subs}) give the SC factorization of
(\ref{rk}) as%
\begin{align*}
t_{n+1} &  =(1)\left(  \frac{a_{n}(t_{n}t_{n-1}\ldots t_{n-k+2})^{-1}}%
{(t_{n}t_{n-1}\ldots t_{n-k+1})^{-1}}+b_{n}\right)  =a_{n}t_{n-k+1}+b_{n},\\
x_{n+1} &  =t_{n+1}x_{n}.
\end{align*}

In this case, the factor equation is linear non-homogeneous with a time delay
of $k-1$ and can be solved to obtain an explicit solution of (\ref{rk})
through (\ref{3}), if desired. Alternatively, we can quickly derive
information about the asymptotic behavior of (\ref{rk}). For instance, if%
\[
\lim_{n\rightarrow\infty}a_{n}=a>0,\ \lim_{n\rightarrow\infty}b_{n}%
=b\geq0,\ a+b\neq1
\]

\noindent then we conclude that all positive solutions of (\ref{rk}) converge
to zero if $a+b<1$ and to $\infty$ if $a+b>1.$

\medskip

\noindent\textbf{Example 2}. This example illustrates a situation where
(\ref{kth}) and (\ref{1aa}) are both HD1, although with respect to different
groups. Consider the third order equation%
\begin{equation}
x_{n+1}=x_{n}+\frac{a(x_{n}-x_{n-1})^{2}}{x_{n-1}-x_{n-2}},\quad
a\neq0.\label{2hd1}%
\end{equation}

Relative to the additive group $G=\mathbb{R}$, this equation is HD1 with the
evident form $x_{n}-x_{n-1}$. Specifically, (\ref{2hd1}) has the form
symmetry
\[
H(u,v,w)=[u-v,v-w],\quad h(t)=-t.
\]

Making the substituion $t_{n}=x_{n}-x_{n-1}$, or using (\ref{subsum}) we get
the SC factorization%
\begin{align}
t_{n+1}  &  =\frac{at_{n}^{2}}{t_{n-1}},\quad t_{0}=x_{0}-x_{-1}%
,\ t_{-1}=x_{-1}-x_{-2}\label{hd1m}\\
x_{n+1}  &  =t_{n+1}+x_{n}.\nonumber
\end{align}

This is a type-(2,1) order reduction. Note that $t_{n}\neq0$ for $n\geq-1$ if
initial values satisfy
\begin{equation}
x_{0},x_{-2}\neq x_{-1}. \label{iv}%
\end{equation}

Relative to the multiplicative group of all nonzero real numbers, the second
order equation (\ref{hd1m}) is HD1 with form symmetry%
\[
H(u,v)=\frac{u}{v},\quad h(s)=\frac{1}{s}.
\]

Making the substitution (\ref{subs}) gives the type-(1,1) order reduction
\begin{align*}
s_{n+1}  &  =as_{n},\quad s_{0}=\frac{t_{0}}{t_{-1}}\\
t_{n+1}  &  =s_{n+1}t_{n}%
\end{align*}

Now using (\ref{3}) and (\ref{3a}) we obtain the following formula for
solutions of (\ref{2hd1}) subject to (\ref{iv}):
\[
x_{n}=x_{0}+t_{0}\sum_{j=1}^{n}s_{0}^{j}a\,^{j(j+1)/2},\quad s_{0}=\frac
{t_{0}}{t_{-1}}=\frac{x_{0}-x_{-1}}{x_{-1}-x_{-2}}.
\]

\medskip

\noindent\textbf{Example 3}. Consider the following variant of Eq.(\ref{2hd1}%
):%
\begin{equation}
x_{n+1}=x_{n}+\frac{a(x_{n}-x_{n-1})}{x_{n-1}-x_{n-2}},\quad a\neq0\label{hs}%
\end{equation}
subject to (\ref{iv}). As in Example 2, the HD1 form $x_{n}-x_{n-1}$ gives the
SC factorization%
\begin{align}
t_{n+1} &  =\frac{at_{n}}{t_{n-1}},\quad t_{0}=x_{0}-x_{-1},\ t_{-1}%
=x_{-1}-x_{-2}\label{hd0}\\
x_{n+1} &  =t_{n+1}+x_{n}.\label{hd0a}%
\end{align}

Unlike (\ref{hd1m}), Eq.(\ref{hd0}) is not HD1. But a straightforward
calculation shows that every solution of (\ref{hd0}) has period 6 as follows%
\begin{equation}
\left\{  t_{-1},t_{0},\frac{at_{0}}{t_{-1}},\frac{a^{2}}{t_{-1}},\frac{a^{2}%
}{t_{0}},\frac{at_{-1}}{t_{0}}\right\}  .\label{p6}%
\end{equation}

Thus we may use (\ref{hd0a}) and (\ref{3a}) to calculate the corresponding
solution of (\ref{hs}) explicitly: For each $n,$ there are integers
$\delta_{n}\geq0$ and $0\leq\rho_{n}\leq5$ such that $n=6\delta_{n}+\rho_{n}.$
If $\sigma$ is the sum of the six numbers in (\ref{p6}) then the explicit
solution of Eq.(\ref{hs}) may be stated as%
\[
x_{n}=x_{0}+%
%TCIMACRO{\tsum _{i=1}^{n}}%
%BeginExpansion
{\textstyle\sum_{i=1}^{n}}
%EndExpansion
t_{i}=x_{0}+\frac{\sigma}{6}(n-\rho_{n})+%
%TCIMACRO{\tsum _{i=n-\rho_{n}+1}^{n}}%
%BeginExpansion
{\textstyle\sum_{i=n-\rho_{n}+1}^{n}}
%EndExpansion
t_{i}%
\]
where every $t_{i}$ is in the set (\ref{p6}) and the last term is zero if
$\rho_{n}=0.$

\pagebreak

\noindent\textbf{Remark}. (\textit{The full triangular factorization property})

The equation in Example 2 has an interesting extra feature: it can be fully SC
factored as a system of first order difference equations%
\begin{align*}
s_{n+1}  &  =as_{n}\\
t_{n+1}  &  =s_{n+1}t_{n}=as_{n}t_{n}\\
x_{n+1}  &  =t_{n+1}+x_{n}=as_{n}t_{n}+x_{n}%
\end{align*}

This system is \textit{triangular} in the sense that each equation is
independent of the variables in the equations below it. For a general
discussion of the periodic solutions of systems of this type see \cite{al},
\cite{kloe}.

If a difference equation of order $k+1$ has the property that it can be
factored completely into a triangular system of first order difference
equations then we say that the difference equation has the \textit{full
triangular factorization} \textit{property} or that it is FTF. Clearly, every
HD1 equation of order 2 is FTF but it is by no means clear if all HD1
equations of order 3 or greater are FTF. For instance, it is not obvious that
Eq.(\ref{hs}) in Example 3 does in fact have the FTF property. The difficulty
there is due to the non-HD1 nature of (\ref{hd0}) which leads to a form
symmetry that involves complex functions (see Example 5 below).

For non-HD1 equations, the FTF property is not clear even for equations of
order 2. But in Corollary 1 in the next section we show that every
\textit{linear} non-homogeneous equation of order $k+1$ is FTF with a complete
factorization into a triangular system of linear non-homogeneous first order equations.

\section{Separability and type-$(1,k)$ factorization}

The class of HD1 functions does not include certain familiar functions. For
example, the linear non-homogeneous function $\phi(u,v)=au+bv+c$ is HD1
relative to the group of all real numbers under addition only when $a+b=1$; it
is HD1 relative to the group $(0,\infty)$ under ordinary multiplication only
when $c=0$ and $a,b\geq0$. These restrictions suggest that a proper study of
order reducible form symmetries for linear difference equations does not
belong in the context of HD1 equations.

In this section we define a class of equations that properly includes all
linear non-homogeneous difference equations with constant coefficients as well
as some other interesting non-HD1 equations. Before discussing this class,
recall that a type-$(1,k)$ equation has a SC factorization with factor of
order $m=1$ and cofactor of order $k.$ Therefore, the form symmetry is a
scalar function that may be written as
\[
H(u_{0},\ldots,u_{k})=u_{0}\ast h(u_{1},\ldots,u_{k}).
\]

Now we define a function $\phi:G^{k+1}\rightarrow G$ to be \textit{separable}
(or \textit{algebraically factorable)} relative to $G$ if there are $k+1$
functions $\phi_{j}:G\rightarrow G$, $j=0,1,\ldots,k$ such that for all
$u_{0},\ldots,u_{k}\in G,$
\[
\phi(u_{0},\ldots,u_{k})=\phi_{0}(u_{0})\ast\cdots\ast\phi_{k}(u_{k}).
\]

Note that every linear non-homogeneous function is trivially separable
relative to every additive subgroup of the complex numbers $\mathbb{C}$. The
rational function $\phi(u,v)=au^{p}/v$ is separable relative to the group of
nonzero real numbers under multiplication for every integer $p$ but it is HD1
relative to the same group if and only if $p=2.$ The exponential function
$\phi(u,v)=ve^{a-bu-cv}$ is separable relative to the group of non-zero real
numbers under multiplication but it is not HD1 relative to that group.

\subsection{Additive forms}

We define Eq.(\ref{kth}) to be separable if every function $f_{n}$ is
separable relative to the underlying group $G.$ In this section we consider
the following separable version of (\ref{kth}) over the group of complex
numbers $\mathbb{C}$ under addition:%
\begin{equation}
x_{n+1}=\alpha_{n}+\phi_{0}(x_{n})+\phi_{1}(x_{n-1})+\cdots\phi_{k}(x_{n-k}).
\label{mde}%
\end{equation}
with%
\begin{equation}
x_{-j},\alpha_{n}\in\mathbb{C},\quad\phi_{j}:\mathbb{C}\rightarrow
\mathbb{C},\quad j=0,1,\ldots,k. \label{mdec}%
\end{equation}

It is not strictly necessary for the sake of applications that the maps
$\phi_{j}$ be defined on all of $\mathbb{C}$\ (indeed, in most applications
they are defined on the set $\mathbb{R}$ of real numbers) but we make a strong
assumption to reduce the amount of technical details in this article. The use
of complex numbers is necessary because form symmetries of (\ref{mde}) may be
complex even if all quantities in (\ref{mdec}) are real (this happens in
particular for linear equations).

The next result from \cite{fs2} shows that Eq.(\ref{mde}) has an
order-reducing form symmetry if one of the $k+1$ functions $\phi_{0}%
,\ldots,\phi_{k}$ can be expressed as a particular linear combination of the
remaining $k$ functions. The form symmetry in this case gives a type-$(1,k)$
order reduction of (\ref{mde}).

\medskip

\noindent\textbf{Theorem 3}. \textit{Assume that there is a constant }%
$c\in\mathbb{C}$\textit{ such that the functions }$\phi_{0},\ldots,\phi_{k}%
$\textit{ in Eq.(\ref{mde}) satisfy}%

\begin{equation}
c^{k+1}z-c^{k}\phi_{0}(z)-c^{k-1}\phi_{1}(z)-\cdots-c\phi_{k-1}(z)-\phi
_{k}(z)=0\quad\text{for all }z.\label{lc}%
\end{equation}

\textit{Then (\ref{mde}) has the following form symmetry}%

\begin{equation}
H(z_{0},z_{1},\ldots,z_{k})=z_{0}+h_{1}(z_{1})+\cdots+h_{k}(z_{k})
\label{hsep}%
\end{equation}
\textit{where }%

\begin{equation}
h_{j}(z)=c^{j}z-c^{j-1}\phi_{0}(z)-\cdots-\phi_{j-1}(z),\quad j=1,\ldots k
\label{hjz}%
\end{equation}

\textit{The form symmetry in (\ref{hsep}) and (\ref{hjz}) yields the
type-}$(1,k)$\textit{ order reduction}%

\begin{align}
z_{n+1}  &  =\alpha_{n}+cz_{n},\quad z_{0}=x_{0}+h_{1}(x_{-1})+\cdots
+h_{k}(x_{-k})\label{lnh}\\
x_{n+1}  &  =z_{n+1}-h_{1}(x_{n})-\cdots-h_{k}(x_{n-k+1}). \label{ro}%
\end{align}

\medskip

\noindent\textbf{Remark. }Note that the factor equation (\ref{lnh}) has order
1 and the cofactor (\ref{ro}) has order $k$ in this case. For reference, we
note that (\ref{lc}) and (\ref{hjz}) imply the following%
\begin{equation}
ch_{k}(z)=\phi_{k}(z). \label{hk}%
\end{equation}

\medskip

A significant feature of Eq.(\ref{ro}) is that it has the same form as
(\ref{mde}). Thus if the functions $h_{1},\ldots,h_{k}$ satisfy the analog of
(\ref{lc}) for some constant $c^{\prime}\in\mathbb{C}$ then Theorem 3 can be
applied to (\ref{ro}). The next result exploits this feature by applying
Theorem 3 to a linear non-homogeneous equation repeatedly until we are left
with a triangular system of first order linear equations.

\medskip

\noindent\textbf{Corollary 1}. \textit{The linear non-homogeneous difference
equation of order }$k+1$\textit{ with constant coefficients}%
\begin{equation}
x_{n+1}+b_{0}x_{n}+b_{1}x_{n-1}+\cdots+b_{k}x_{n-k}=\alpha_{n}\label{LDE}%
\end{equation}
\textit{where }$b_{0},\ldots,b_{k}$\textit{, }$\alpha_{n}\in\mathbb{C}%
$\textit{ has the FTF property and is equivalent to the following triangular
system of }$k+1$\textit{ first order linear non-homogeneous equations}%
\begin{align*}
z_{0,n+1} &  =\alpha_{n}+c_{0}z_{0,n},\\
z_{1,n+1} &  =z_{0,n+1}+c_{1}z_{1,n}\\
&  \vdots\\
z_{k,n+1} &  =z_{k-1,n+1}+c_{k}z_{k,n}%
\end{align*}
\textit{in which }$z_{k,n}=x_{n}$\textit{ is the solution of Eq.(\ref{LDE})
and the constants }$c_{0},c_{1},\ldots,c_{k}$\textit{ are the eigenvalues of
the homogeneous part of (\ref{LDE}), i.e., roots of the characteristic
polynomial}%
\begin{equation}
P(z)=z^{k+1}+b_{0}z^{k}+b_{1}z^{k-1}+\cdots+b_{k-1}z+b_{k}.\label{cp}%
\end{equation}

\textbf{Proof}. Defining $\phi_{j}(z)=-b_{j}z$ for $j=1,\ldots k$ and applying
Theorem 3 above yields the SC factorization%
\begin{align}
z_{0,n+1}  &  =\alpha_{n}+c_{0}z_{0,n}\nonumber\\
x_{n+1}  &  =z_{0,n+1}-\beta_{1,0}x_{n}-\cdots-\beta_{1,k-1}x_{n-k+1}
\label{l1}%
\end{align}
where $c_{0}$ satisfies (\ref{lc})%
\[
c_{0}^{k+1}z+c_{0}^{k}b_{0}z+c_{0}^{k-1}b_{1}z+\cdots+c_{0}b_{k-1}z+b_{k}z=0
\]
for all $z\in\mathbb{C},$ i.e. $c_{0}$ is a root of the characteristic
polynomial $P$ in (\ref{cp}). Further, the numbers $\beta_{1,j}$ are given via
the function $h_{j}$ in (\ref{hjz}) and (\ref{hk}) as%
\[
h_{j}(z)=\beta_{1,j-1}z,\quad\beta_{1,j-1}=c_{0}^{j}+c_{0}^{j-1}b_{0}%
+\cdots+b_{j-1},\ c_{0}\beta_{1,k-1}=-b_{k-1.}%
\]

Alternatively, the numbers $\beta_{1,j}$ may be calculated from the recursion%
\begin{equation}
\beta_{1,j}=c_{0}\beta_{1,j-1}+b_{j},\quad j=1,\ldots k-1,\ \beta_{1,0}%
=c_{0}+b_{0},\ c_{0}\beta_{1,k-1}=-b_{k-1}. \label{b1j}%
\end{equation}

Next, since Eq.(\ref{l1}), i.e., $x_{n+1}+\beta_{1,0}x_{n}+\cdots
+\beta_{1,k-1}x_{n-k+1}=z_{0,n+1}$ is of the same type as (\ref{LDE}), Theorem
3 can be applied to it to yield the SC factorization%
\begin{align*}
z_{1,n+1}  &  =z_{0,n+1}+c_{1}z_{1,n}\\
x_{n+1}  &  =z_{1,n+1}-\beta_{2,0}x_{n}-\cdots-\beta_{2,k-2}x_{n-k+2}%
\end{align*}
in which $c_{1}$ satisfies (\ref{lc}) for (\ref{l1}), i.e., the power is
reduced by 1 and coefficients adjusted appropriately as in%
\[
c_{1}^{k}+\beta_{1,0}c_{1}^{k-1}+\beta_{1,1}c_{1}^{k-2}+\cdots+\beta
_{1,k-2}c_{0}+\beta_{1,k-1}=0.
\]

Now we show that $c_{1}$ is also a root of $P$ in (\ref{cp}). Define%
\[
P_{1}(z)=z^{k}+\beta_{1,0}z^{k-1}+\cdots+\beta_{1,k-2}z+\beta_{1,k-1}%
\]
so that $c_{1}$ is a root of $P_{1}.$ If it is shown that $(z-c_{0}%
)P_{1}(z)=P(z)$ then $P_{1}$ divides $P$ so $c_{1}$ is a root of $P.$ Direct
calculation using (\ref{b1j}) shows%
\begin{align*}
(z-c_{0})P_{1}(z)  &  =z^{k+1}+\beta_{1,0}z^{k}+\beta_{1,1}z^{k}+\cdots
+\beta_{1,k-1}z\\
&  -c_{0}z^{k}-c_{0}\beta_{1,0}z^{k-1}-\cdots-c_{0}\beta_{1,k-2}z-c_{0}%
\beta_{1,k-1}\\
&  =z^{k+1}+(c_{0}+b_{0})z^{k}+(c_{0}\beta_{1,0}+b_{1})z^{k-1}\cdots
+(c_{0}\beta_{1,k-2}+b_{k-2})z\\
&  -c_{0}z^{k}-c_{0}\beta_{1,0}z^{k-1}-\cdots-c_{0}\beta_{1,k-2}z-c_{0}%
\beta_{1,k-1}\\
&  =P(z).
\end{align*}

Therefore, the above process inductively generates the system in the statement
of this corollary.

\medskip

\noindent\textbf{Remarks}. \textit{(Operator factorization, complementary and
particular solutions)}

1. The triangular SC factorization of Corollary 1 is essentially what is
obtained through operator factorization. If $Ex_{n}=x_{n+1}$ represents the
forward shift operator then as is well-known, the eigenvalues factor the
operator $P(E)$ with $P$ defined by (\ref{cp}); i.e., (\ref{LDE}) can be
written as%
\begin{equation}
(E-c_{0})(E-c_{1})\cdots(E-c_{k})x_{n-k}=\alpha_{n}. \label{op}%
\end{equation}

Now if we define
\begin{equation}
(E-c_{1})\cdots(E-c_{k})x_{n-k}=y_{0,n}\label{op1}%
\end{equation}
then (\ref{op}) can be written as%
\[
y_{0,n+1}-c_{0}y_{0,n}=\alpha_{n}%
\]
which is the first equation in the triangular system of Corollary 1 with
$y_{0,n}=z_{0,n}$. We may continue in this fashion by applying the same idea
to (\ref{op1}); we set%
\[
(E-c_{2})\cdots(E-c_{k})x_{n-k}=y_{1,n}%
\]
and write (\ref{op1}) as $y_{1,n+1}-c_{1}y_{1,n}=z_{0,n}$ which is the second
equation in the triangular system if $y_{1,n}=z_{1,n-1}.$ The reduction in the
time index $n$ here is due to the removal of one occurrence of $E.$ Proceeding
in this fashion, setting $y_{j,n}=z_{j,n-j}$ at each step, we eventually
arrive at
\[
(E-c_{k})x_{n-k}=y_{k-1,n}\Rightarrow x_{n+1-k}=y_{k-1,n}+c_{k}x_{n-k}.
\]

Thus, with $y_{k-1,n}=z_{k-1,n-k+1}$ the preceding equation is the same as the
last equation in the system of Corollary 1.

2. With Corollary 1 we may obtain the eigenvalues and both the particular
solution and the solution of the homogeneous part of (\ref{LDE})
simultaneously without needing to guess linearly independent solutions,
namely, the complex exponentials. We indicate how this is done in the second
order case $k=1$ which is also representative of the higher order cases.
First, for a given sequence $s=\{s_{n}\}$ of complex numbers and for each
$c\in\mathbb{C}$, let us define the quantity%
\[
\sigma_{n}(s;c)=\sum_{j=1}^{n}c^{j-1}s_{n-j}%
\]
and note that for sequences $s,t$ and numbers $a,b\in\mathbb{C}$, $\sigma
_{n}(as+bt;c)=a\sigma_{n}(s;c)+b\sigma_{n}(t;c),$ i.e., $\sigma_{n}(\cdot,c)$
is a linear operator on the space of complex sequences for each $n\geq1$ and
each $c\in\mathbb{C}.$ Further, if $s_{n}=ab^{n}$ then it is easy to see that%
\begin{equation}
\sigma_{n}(s;c)=\left\{
\begin{array}
[c]{ll}%
a(b^{n}-c^{n})/(b-c), & c\neq b\\
anb^{n-1}, & c=b
\end{array}
\right.  .\label{sig}%
\end{equation}

Now, if $k=1$ then the semiconjugate factorization of (\ref{LDE}) into first
order equations is%

\begin{align}
z_{n+1}  &  =\alpha_{n}+c_{0}z_{n},\quad z_{0}=x_{0}+(c_{0}-b_{0}%
)x_{-1}\label{l2a}\\
x_{n+1}  &  =z_{n+1}+c_{1}x_{n}. \label{l2b}%
\end{align}

A straightforward inductive argument gives the solution of (\ref{l2a}) as%
\begin{equation}
z_{n}=z_{0}c_{0}^{n}+\sigma_{n}(\alpha;c_{0}). \label{sl2a}%
\end{equation}

Next, insert (\ref{sl2a}) into (\ref{l2b}), set $\gamma_{n}=z_{n+1}$ and
repeat the above argument to obtain the general solution of (\ref{LDE}) for
$k=1$, i.e., $x_{n}=x_{0}c_{1}^{n}+\sigma_{n}(\gamma;c_{1}).$ If $c_{1}\neq
c_{0}$ then from (\ref{sig}) we obtain after combining some terms and noting
that $\gamma_{0}=z_{1}=\alpha_{0}+c_{0}z_{0}$,
\[
x_{n}=\left(  \frac{\alpha_{0}+c_{0}z_{0}}{c_{0}-c_{1}}\right)  c_{0}%
^{n}+\left(  x_{0}-\frac{\alpha_{0}+c_{0}z_{0}}{c_{0}-c_{1}}\right)  c_{1}%
^{n}+\sigma_{n}(\sigma^{\prime}(\alpha;c_{0});c_{1}).
\]
where $\sigma^{\prime}=\{\sigma_{n+1}(\alpha;c_{0})\}.$ We recognize the first
two terms of the above sum as giving the solution of the homogeneous part of
(\ref{LDE}) and the last term as giving the particular solution. In the case
of repeat eigenvalues, i.e., $c_{1}=c_{0}$ again from (\ref{sig}) we get
\[
x_{n}=[x_{0}c_{0}+(\alpha_{0}+c_{0}z_{0})n]c_{0}^{n-1}+\sigma_{n}%
(\sigma^{\prime}(\alpha;c_{0});c_{0}).
\]

\subsection{Multiplicative forms}

As another application of Theorem 3 we consider the following difference
equation on the positive real line%
\begin{gather}
y_{n+1}=\beta_{n}\psi_{0}(y_{n})\psi_{1}(y_{n-1})\cdots\psi_{k}(y_{n-k}%
),\label{mult}\\
\beta_{n},y_{-j}\in(0,\infty),\ \psi_{j}:(0,\infty)\rightarrow(0,\infty
),\ j=0,\ldots k.\nonumber
\end{gather}

Taking the logarthim of Eq.(\ref{mult}) changes its multiplicative form to an
additive one. Specifically by defining%
\[
x_{n}=\ln y_{n},\ y_{n}=e^{x_{n}},\ \ln\beta_{n}=\alpha_{n},\ \phi_{j}%
(r)=\ln\psi_{j}(e^{r}),\ j=0,\ldots k,\ r\in\mathbb{R}%
\]
we can transform (\ref{mult}) into (\ref{mde}). Then Theorem 3 implies the
following generalization of the main result of \cite{fs1}.

\medskip

\noindent\textbf{Corollary 2}. \textit{Eq.(\ref{mult}) has a form symmetry }%
\begin{equation}
H(t_{0},t_{1},\ldots,t_{k})=t_{0}h_{1}(t_{1})\cdots h_{k}(t_{k}) \label{Hmult}%
\end{equation}
\textit{if there is }$c\in\mathbb{C}$\textit{ such that the following is true
for all }$t>0,$%
\begin{equation}
\psi_{0}(t)^{c^{k}}\psi_{1}(t)^{c^{k-1}}\cdots\psi_{k}(t)=t^{c^{k+1}}.
\label{lcmult}%
\end{equation}

\textit{The functions }$h_{j}$ \textit{in (\ref{Hmult}) are given as}%
\begin{equation}
h_{j}(t)=t^{c^{j}}\psi_{0}(t)^{-c^{j-1}}\cdots\psi_{j-1}(t),\quad j=1,\ldots k
\label{hjmult}%
\end{equation}
\textit{and the form symmetry in (\ref{Hmult}) and (\ref{hjmult}) yields the
type-}$(1,k)$\textit{ order reduction}
\begin{subequations}
\label{or}%
\begin{align}
r_{n+1}  &  =\beta_{n}r_{n}^{c},\quad r_{0}=y_{0}h_{1}(y_{-1})\cdots
h_{k}(y_{-k})\label{fmult}\\
y_{n+1}  &  =\frac{r_{n+1}}{h_{1}(y_{n})\cdots h_{k}(y_{n-k+1})}.
\label{flink}%
\end{align}

\medskip

\noindent\textbf{Example 4 (}A simple equation with complicated multistable
solutions). Equations of type (\ref{mde}) or (\ref{mult}) are capable of
exhibiting complex behavior, including the generation of coexisting
\textit{stable} solutions of many different types that range from periodic to
chaotic. As a specific example of such multistable equations consider the
following second-order equation%
\end{subequations}
\begin{equation}
x_{n+1}=x_{n-1}e^{a-x_{n}-x_{n-1}},\quad x_{-1},x_{0}>0.\label{exp}%
\end{equation}

Note that Eq.(\ref{exp}) has up to two isolated fixed points. One is the
origin which is repelling if $a>0$ (eigenvalues of linearization are $\pm
e^{a/2}$) and the other fixed point is $\bar{x}=a/2$. If $a>4$ then $\bar{x}$
is unstable and non-hyperbolic because the eigenvalues of the linearization of
(\ref{exp}) are $-1$ and $1-a/2$. The computer-generated diagram in Figure 1
shows the variety of stable periodic and non-periodic solutions that occur
with $a=4.6$ and one initial value $x_{-1}=2.3$ fixed and the other initial
value $x_{0}$ changing from 2.3 to 4.8; i.e., approaching (or moving away
from) the fixed point $\bar{x}$ on a straight line segment in the plane.

%************************* FIGURE 1 *****************************
\begin{figure}[ptb]
\begin{center}
\includegraphics[width=4.795in]{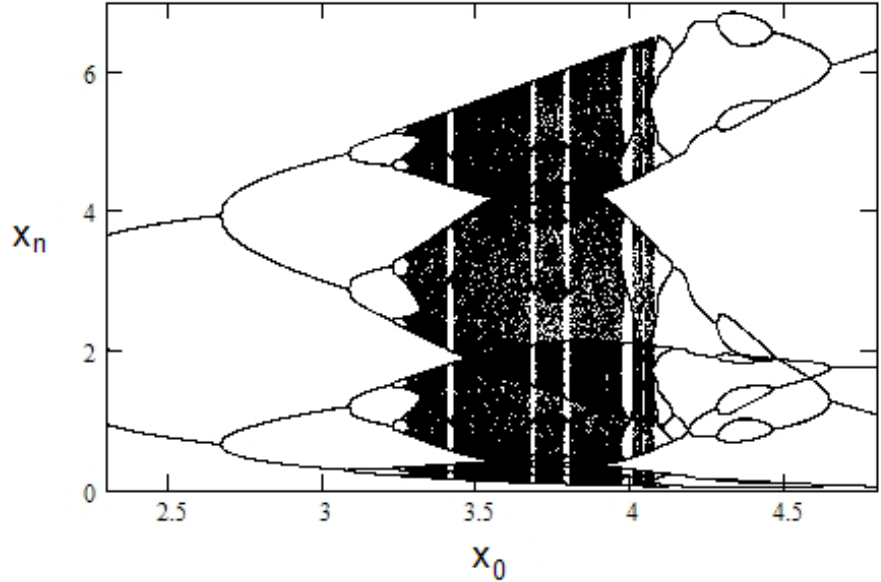}
\end{center}
\caption{{\small {Bifurcations of solutions of Eq.(\ref{exp}) with a changing
initial value; $a=4.6$ is fixed.}}}%
\end{figure}

In Figure 1, for every grid value of $x_{0}$ in the range 2.3-4.8, the last
200 (of 300) points of the solution $\{x_{n}\}$ are plotted vertically. In
this figure, stable solutions with periods 2, 4, 8, 12 and 16 can be easily
identified. All of the solutions that appear in Figure 1 represent
\textit{coexisting stable orbits} of Eq.(\ref{exp}). There are also periodic
and non-periodic solutions which do not appear in Figure 1 because they are
unstable (e.g., the fixed point $\bar{x}=2.3$). Additional bifurcations of
both periodic and non-periodic solutions occur outside the range 2.3-4.8 which
are not shown in Figure 1.

Understanding the behavior for solutions of Eq.(\ref{exp}) is made easier when
we look at its SC factorization given by (\ref{fmult}) and (\ref{flink}). Here
$k=1$ and
\[
\psi_{0}(t)=e^{-t},\quad\psi_{1}(t)=te^{-t},\quad\beta_{n}=e^{a}\text{ for all
}n.
\]

Thus (\ref{lcmult}) takes the form%
\begin{align*}
\psi_{0}(t)^{c}\psi_{1}(t)  &  =t^{c^{2}}\text{ for all }t>0\\
e^{-ct}te^{-t}  &  =t^{c^{2}}\text{ for all }t>0
\end{align*}

The last equality is true if $c=-1,$ which leads to the form symmetry%
\[
h_{1}(t)=t^{-1}\psi_{0}(t)^{-1}=\frac{1}{te^{-t}}\Rightarrow H(u_{0}%
,u_{1})=\frac{u_{0}}{u_{1}e^{-u_{1}}}%
\]
and SC factorization%
\begin{align}
r_{n+1} &  =\frac{e^{a}}{r_{n}},\quad r_{0}=x_{0}h_{1}(x_{-1})=\frac{x_{0}%
}{x_{-1}e^{-x_{-1}}}\label{star1}\\
x_{n+1} &  =\frac{r_{n+1}}{h_{1}(x_{n})}=r_{n+1}x_{n}e^{-x_{n}}.\label{star2}%
\end{align}

All solutions of (\ref{star1}) with $r_{0}\neq e^{a/2}$ are periodic with
period 2:%
\[
\left\{  r_{0},\frac{e^{a}}{r_{0}}\right\}  =\left\{  \frac{x_{0}}%
{x_{-1}e^{-x_{-1}}},\frac{x_{-1}e^{a-x_{-1}}}{x_{0}}\right\}  .
\]

Hence the orbit of each nontrivial solution $\{x_{n}\}$ of (\ref{exp}) in the
plane is restricted to the pair of curves%
\begin{equation}
\xi_{1}(t)=\frac{e^{a}}{r_{0}}te^{-t}\quad\text{and\quad}\xi_{2}%
(t)=r_{0}te^{-t}.\label{ic}%
\end{equation}

Now, if $x_{-1}$ is fixed and $x_{0}$ changes, then $r_{0}$ changes
proportionately to $x_{0}$. These changes in initial values are reflected as
changes in \textit{parameters} in (\ref{star2}). The orbits of the one
dimensional map $bte^{-t}$ where $b=r_{0}$ or $e^{a}/r_{0}$\ exhibit a variety
of behaviors as the parameter $b$ changes according to well-known rules such
as the fundamental ordering of cycles and the occurrence of chaotic behavior
with the appearnce of period-3 orbits when $b$ is large enough; see, e.g.,
\cite{bc}, \cite{ce}, \cite{li}, \cite{hsbk}. Eq.(\ref{star2}) splits these
behaviors evenly over the pair of curves (\ref{ic}) as the initial value
$x_{0}$ changes; see Figure 2 which shows the orbits of (\ref{exp}) for two
different initial values $x_{0}$ with $a=4.6;$ the first 100 points of each
orbit are discarded in these images so as to highlight the asymptotic behavior
of each orbit. The splitting over the pair of curves $\xi_{1},\xi_{2}$ also
explains why odd periods do not appear in Figure 1.

%********************* FIGURE 2 ************************************
\begin{figure}[ptb]
\begin{center}
\includegraphics[width=4.557in]{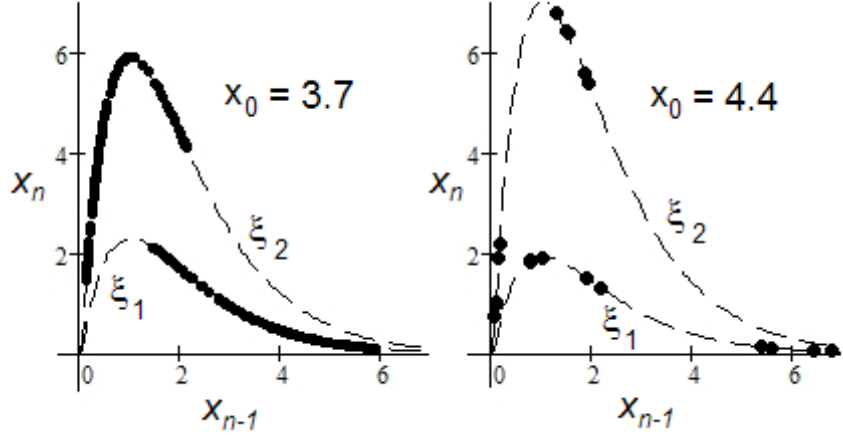}
\end{center}
\caption{{\small {Two of the orbits in Figure 1 shown here on their loci of
two curves $\xi_{1}$, $\xi_{2}$ in the state space.} }}%
\end{figure}

\medskip

\noindent\textbf{Example 5}. We now combine different types of form symmetry
to show that Eq.(\ref{hs}) in Example 3 has the FTF property. If we assume
that
\begin{equation}
a>0,\ x_{0}>x_{-1}>x_{-2},\ x_{0}\geq0\label{prest}%
\end{equation}
then the multiplicative group $(0,\infty)$ of positive real numbers is the
invariant set of Eq.(\ref{hd0}). Since (\ref{hd0}) is obviously separable over
$(0,\infty)$, we check equality (\ref{lcmult}) in Corollary 2 with $\psi
_{0}(t)=t$ and $\psi_{1}(t)=1/t$:%
\[
t^{c-1}=t^{c^{2}}\quad\text{for all }t>0.
\]

This condition holds if $c^{2}-c+1=0.$ The quadratic has complex roots%
\[
c_{\pm}=\frac{1\pm i\sqrt{3}}{2}%
\]
so by Corollary 2, Eq.(\ref{hd0}) has a form symmetry%
\[
H(u_{0},u_{1})=u_{0}h_{1}(u_{1}),\quad\text{where }h_{1}(t)=t^{c_{+}}%
t^{-1}=t^{-c_{-}}%
\]
and an SC factorization%
\begin{align}
r_{n+1} &  =ar_{n}^{c_{+}},\quad r_{0}=t_{0}h_{1}(t_{-1})\label{3r}\\
t_{n+1} &  =r_{n+1}t_{n}^{c_{-}}.\label{3t}%
\end{align}

The three equations (\ref{3r}), (\ref{3t}) and (\ref{hd0a}) establish that
Eq.(\ref{hs}) has the FTF property with a factorization%
\begin{align*}
r_{n+1} &  =ar_{n}^{c_{+}},\\
t_{n+1} &  =ar_{n}^{c_{+}}t_{n}^{c_{-}},\\
x_{n+1} &  =ar_{n}^{c_{+}}t_{n}^{c_{-}}+x_{n}.
\end{align*}

It is noteworthy that in spite of the occurrence of complex exponents, this
system generates positive solutions from positive initial values. This fact
may seem less surprising if Eq.(\ref{hd0}) is transformed into a linear
equation (with complex eigenvalues) by taking logarithms as in the beginning
of this section.

\medskip

\noindent\textbf{Remark}. Under the added restrictions (\ref{prest}),
Corollary 2 may also be applied to Eq.(\ref{hd1m}) of Example 2 to obtain an
SC factorization using the separable type of form symmetry. Is this SC
factorization different from that in Example 2? To see that they are in fact
the same, let $\psi_{0}(t)=t^{2}$ and $\psi_{1}(t)=1/t$ in the equality
(\ref{lcmult}) and require that%
\[
t^{2c-1}=t^{c^{2}}\quad\text{for all }t>0.
\]

This holds if $c$ is a root of the quadratic $c^{2}-2c+1=0,$ i.e., $c=1.$ Thus
using (\ref{hjmult}) we calculate the form symmetry as%
\[
H(u_{0},u_{1})=u_{0}h_{1}(u_{1}),\quad\text{where }h_{1}(t)=t^{1}t^{-2}%
=t^{-1}.
\]

This is just the HD1 form symmetry giving the same SC factorization as in
Example 2.

\vspace{2ex}

\noindent Department of Mathematics, Virginia Commonwealth University,
Richmond, Virginia 23284-2014, USA

\noindent Email: hsedagha@vcu.edu

\end{document}